\documentclass{svjour2}                    
\smartqed  
\usepackage{mathptmx}      
\newcommand{\mathbb}[1]{#1}
\newcommand{\be}{\begin{equation}}
\newcommand{\ee}{\end{equation}}
\newcommand{\bea}{\begin{eqnarray}}
\newcommand{\eea}{\end{eqnarray}}
\newcommand{\rd}[1]{{\rm d}#1\,}

\newcommand{\rdf}[3]{\frac{{\rm d}^{#1}#2}{#3}\,}

\newcommand{\re}[1]{{\rm e}^{#1}}
\newcommand{\bv}[1]{#1}

\newcommand{\eref}[1]{Eq.~(\ref{#1})}
\newcommand{\Eref}[1]{Equation (\ref{#1})}
\newcommand{\sref}[1]{Sec.~\ref{#1}}

\begin{document}

\title{Radiative Transport Equation for Bloch Electrons 
in Electromagnetic Fields
}


\author{Manabu Machida}


\institute{M. Machida \at
 Department of Mathematics, University of Michigan,\\
 Ann Arbor, MI 48109, USA \\
 \email{mmachida@umich.edu}           
}

\date{Received: date / Accepted: date}

\maketitle

\begin{abstract}
The radiative transport equation for the Schr\"{o}dinger equation in 
a periodic potential with a weak random potential 
in electromagnetic fields is derived using 
asymptotic expansion.

\end{abstract}

\section{Introduction}
\label{intro}

The radiative transport equation (RTE) is a linear Boltzmann equation, which 
describes propagation of energy density 
in a random medium \cite{Chandrasekhar60,Case67}.  
Rigorous analysis on the transport limit for the Schr\"{o}dinger equation 
with random potential has been studied 
\cite{Martin75,SpohnH77,DellAntonio83,Ho93,Markowich94,Gerard97,Erdos00,Ryzhik96,Bal02,Bal99}.  In particular, the RTE was obtained from the 
Schr\"{o}dinger equation with time independent \cite{Ryzhik96} and 
time dependent \cite{Bal02} random potential.  
In \cite{Bal99}, the RTE was derived from the Schr\"{o}dinger equation with 
random and periodic potential.  See also recent review \cite{Bal10} and 
references therein.  
In the absence of random potential, the semiclassical equations of motion 
for the Schr\"{o}dinger equation with periodic potential and 
electromagnetic fields has been considered 
\cite{Guillot88,Chang96,Sundaram99,Hovermann01,Dimassi02,Panati03}.  

In this paper, we consider noninteracting electrons in a periodic 
potential (Bloch electrons) with a weak random potential and 
apply electromagnetic fields.  For this system, the RTE 
is derived in the same way as \cite{Bal99} except that we here take 
electromagnetic fields into account.  

This paper is organized as follows.  In \sref{sec:schr}, a two-scale 
asymptotic expansion is introduced for the Wigner distribution function.  
In \sref{sec:bloch}, Bloch functions are considered.  
The Wigner distribution function is decomposed in \sref{sec:w0decomp} and 
\sref{sec:w1decomp} making use of the Bloch functions.  Finally, 
we derive the RTE and obtain our main result \eref{main} in \sref{sec:w0}.

\section{The Schr\"{o}dinger Equation}
\label{sec:schr}

We consider Bloch electrons in an electric field $\bv{E}\in\mathbb{R}^3$ and 
a magnetic field $\bv{B}\in\mathbb{R}^3$.  They are given by 
vector potential $\bv{A}\in\mathbb{R}^3$ and 
scalar potential $\varphi\in\mathbb{R}$ as
\be
\bv{E}=-\nabla\varphi-\frac{\partial}{\partial t}\bv{A},\qquad
\bv{B}=\nabla\times\bv{A}.
\label{vectorpotential}
\ee
The state $\psi(t,\bv{x})\in\mathbb{R}$ ($t\in\mathbb{R}$, 
$\bv{x}\in\mathbb{R}^3$) of this system evolves according to 
the Schr\"{o}dinger equation
\be
i\hbar\frac{\partial}{\partial t}\psi(t,\bv{x})=\mathcal{H}\psi(t,\bv{x}).
\label{schr}
\ee
The Hamiltonian $\mathcal{H}$ is given as
\be
\mathcal{H}=
\frac{1}{2m_{\rm e}}\left(\bv{p}+e\bv{A}\right)^2+e\varphi+
U_{\rm tot}(\bv{x}),
\ee
where $\bv{p}=-i\hbar\nabla_{\bv{x}}$, $-e$ and $m_{\rm e}$ are 
the charge and mass of an electron, and 
$U_{\rm tot}(\bv{x})=U(\bv{x})+V(\bv{x})$.  
Here $U(\bv{x})$ is a $d$-dimensional ($d\le 3$) periodic potential and 
$V(\bv{x})$ is a random potential.  The periodic potential satisfies
\be
U(\bv{z}+\bv{\nu})=U(\bv{z}),
\ee
where $\bv{\nu}$ belongs to the lattice $L$:
\be
L=\left\{\sum_{j=1}^dn_j\bv{e}_j\Biggm|n_j\in\mathbb{Z}\right\}
\ee
and $\bv{e}_1,\dots,\bv{e}_d$ form a basis of $\mathbb{R}^d$ with the dual 
basis $\bv{e}^1,\dots,\bv{e}^d$ defined by
$\bv{e}_j\cdot\bv{e}^k=2\pi\delta_{jk}$.  
The dual lattice $L^*$ is spanned by $\{\bv{e}^k\}$.  
We let $C$ denote the basic period cell of $L$ and $\Omega_{\rm BZ}$ 
denote the Brillouin zone.  Note that
\be
|C|\Omega_{\rm BZ}=(2\pi)^d.
\ee
The random potential has the following properties.
\be
\langle V(\bv{y})V(\bv{y}+\bv{x})\rangle=R(\bv{x}),\quad
\langle\tilde{V}(\bv{q})\tilde{V}(\bv{q}')\rangle=(2\pi)^d
\tilde{R}(\bv{q'})\delta(\bv{q}+\bv{q'}).
\ee
Here $\langle\cdot\rangle$ denotes ensemble average and the Fourier 
transform is defined as
\be
\tilde{V}(\bv{q})=\int_{\mathbb{R}^d}\rd{x}\re{-i\bv{q}\cdot\bv{x}}V(\bv{x}).
\ee
Because of the $U(1)$ gauge symmetry, we can set
\be
\varphi=0.
\ee

We scale the variables as
\be
t\to\frac{t}{\varepsilon},\qquad\bv{x}\to\frac{\bv{x}}{\varepsilon},
\ee
where $\varepsilon$ ($>0$) is small.  We assume weak electromagnetic fields: 
\be
\bv{E}\to\varepsilon\bv{E},\qquad\bv{B}\to\varepsilon\bv{B}.
\ee
\Eref{vectorpotential} implies that $\bv{A}$ is independent of 
$\varepsilon$.  Since the cyclotron radius $\ell$ is inversely proportional 
to $|\bv{B}|$, we have $\ell\to\ell/\varepsilon$ and $\ell$ gets much larger 
than the cell size of the lattice $L$.  Furthermore, 
we assume that the random potential is weak.  We obtain
\bea
\frac{\partial}{\partial t}\psi_{\varepsilon}(t,\bv{x})
&=&
\frac{i\hbar\varepsilon}{2m_{\rm e}}\left[\nabla_{\bv{x}}+
\frac{ie}{\hbar\varepsilon}\bv{A}(t,\bv{x})\right]^2
\psi_{\varepsilon}(t,\bv{x})+
\frac{1}{i\hbar\varepsilon}U\left(\frac{\bv{x}}{\varepsilon}\right)
\psi_{\varepsilon}(t,\bv{x})
\nonumber \\
&+&
\frac{1}{i\hbar\sqrt{\varepsilon}}V\left(\frac{\bv{x}}{\varepsilon}\right)
\psi_{\varepsilon}(t,\bv{x}),
\eea
where the initial wave function 
$\psi_{\varepsilon}(0,\bv{x})\in L^2(R^d)$ is $\varepsilon$-oscillatory.  
In \cite{Bal99}, the RTE is derived from this equation in the absence 
of the vector potential $\bv{A}(t,\bv{x})$.

Let us consider the Wigner distribution function associated with 
$\psi_{\varepsilon}$.
\be
W_{\varepsilon}(t,\bv{x},\bv{k})=
\int_{\mathbb{R}^d}\rdf{}{\bv{y}}{(2\pi)^d}\re{i\bv{k}\cdot\bv{y}}
\psi_{\varepsilon}(t,\bv{x}-\varepsilon\bv{y})
\bar{\psi}_{\varepsilon}(t,\bv{x}).
\label{WDF}
\ee
Whenever necessary, we use the fact that $\bv{k}\in\mathbb{R}^d$ is 
uniquely decomposed as
\be
\bv{k}=\bv{q}+\bv{\mu},\quad
\bv{q}\in\Omega_{\rm BZ},\quad\bv{\mu}\in L^*.
\ee
Note that $W_{\varepsilon}(t,\bv{x},\bv{k})$ and its symmetric version
\be
\int_{\mathbb{R}^d}\rdf{}{\bv{y}}{(2\pi)^d}\re{i\bv{k}\cdot\bv{y}}
\psi_{\varepsilon}\left(t,\bv{x}-\frac{\varepsilon}{2}\bv{y}\right)
\bar{\psi}_{\varepsilon}\left(t,\bv{x}+\frac{\varepsilon}{2}\bv{y}\right)
\ee
have the same weak limit as $\varepsilon\to 0$ \cite{GerardP93}.  
We obtain the time evolution of $W_{\varepsilon}$ as
\bea
\frac{\partial}{\partial t}W_{\varepsilon}(t,\bv{x},\bv{k})+
\Biggl[
\frac{\hbar\bv{k}}{m_{\rm e}}\cdot\nabla_{\bv{x}}+
\frac{i\hbar\varepsilon}{2m_{\rm e}}\Delta_{\bv{x}}+
\frac{e}{m_{\rm e}}\Bigl(\bv{A}(t,\bv{x})\cdot\nabla_{\bv{x}}-
\bv{k}\cdot\bv{A}(t,\nabla_{\bv{\mu}})\Bigr)
\nonumber \\
&\hspace{-16cm}-&\hspace{-8cm}
\frac{e^2}{m_{\rm e}\hbar}
\left(\bv{A}(t,\bv{x})+\frac{i\varepsilon}{2}\bv{A}(t,\nabla_{\bv{\mu}})\right)
\cdot\bv{A}(t,\nabla_{\bv{\mu}})
\Biggr]W_{\varepsilon}(t,\bv{x},\bv{k})
\nonumber \\
&\hspace{-16cm}=&\hspace{-8cm}
\frac{1}{i\hbar\varepsilon}\sum_{\bv{\mu}'\in L^*}
\re{i\bv{\mu}'\cdot\bv{x}/\varepsilon}\tilde{U}(\bv{\mu}')
\left[W_{\varepsilon}(t,\bv{x},\bv{k}-\bv{\mu}')-
W_{\varepsilon}(t,\bv{x},\bv{k})\right]
\nonumber \\
&\hspace{-16cm}+&\hspace{-8cm}
\frac{1}{i\hbar\sqrt{\varepsilon}}\int_{\mathbb{R}^d}
\rdf{}{\bv{k}'}{(2\pi)^d}\re{i\bv{k}'\cdot\bv{x}/\varepsilon}\tilde{V}(\bv{k}')
\left[W_{\varepsilon}(t,\bv{x},\bv{k}-\bv{k}')-
W_{\varepsilon}(t,\bv{x},\bv{k})\right],
\eea
where we used
\be
\psi_{\varepsilon}(t,\bv{x})\frac{\partial}{\partial x_j}
\psi_{\varepsilon}(t,\bv{x})\stackrel{x_j\to\infty}{\rightarrow}0,
\quad(j=1,\dots,d),
\ee
and introduced
\be
U(\bv{y})=\sum_{\bv{\mu}\in L^*}\re{i\bv{\mu}\cdot\bv{y}}\tilde{U}(\bv{\mu}),
\quad
\tilde{U}(\bv{\mu})=
\frac{1}{|C|}\int_C\rd{\bv{y}}\re{-i\bv{\mu}\cdot\bv{y}}U(\bv{y}).
\ee
Note that we have
\be
\frac{1}{|C|}\sum_{\bv{\mu}\in L^*}\re{i\bv{\mu}\cdot\bv{z}}=
\sum_{\bv{\nu}\in L}\delta(\bv{z}-\bv{\nu}).
\label{musum}
\ee

Let us define
\be
H(\bv{x},\bv{\mu},\bv{q})=
\frac{\hbar}{2m_{\rm e}}\left[\bv{k}+\frac{e}{\hbar}\bv{A}(t,\bv{x})\right]^2.
\ee
Noting that
\be
\frac{\partial\bv{A}(\bv{x})}{\partial{x_i}}\frac{\partial}{\partial{\mu_i}}
=
\bv{A}({^t}(0,\dots,0,\frac{\partial}{\partial_{\mu_i}},0,\dots,0)),
\ee
we obtain
\bea
\frac{\partial}{\partial t}W_{\varepsilon}+\left[
\nabla_{\bv{\mu}}H(\bv{x},\bv{\mu},\bv{q})\cdot\nabla_{\bv{x}}-
\nabla_{\bv{x}}H(\bv{x},\bv{\mu},\bv{q})\cdot\nabla_{\bv{\mu}}
\right]W_{\varepsilon}(t,\bv{x},\bv{k})=
\nonumber \\
&&\hspace{-8cm}
\left[-\frac{i\hbar\varepsilon}{2m_{\rm e}}\Delta_{\bv{x}}
+\frac{i\varepsilon e^2}{2m_{\rm e}\hbar}
\bv{A}(t,\nabla_{\bv{\mu}})\cdot\bv{A}(t,\nabla_{\bv{\mu}})\right]
W_{\varepsilon}(t,\bv{x},\bv{k})
\nonumber \\
&\hspace{-16cm}+&\hspace{-8cm}
\frac{1}{i\hbar\varepsilon}\sum_{\bv{\mu}'\in L^*}
\re{i\bv{\mu}'\cdot\bv{x}/\varepsilon}\tilde{U}(\bv{\mu}')
\left[W_{\varepsilon}(t,\bv{x},\bv{k}-\bv{\mu}')-
W_{\varepsilon}(t,\bv{x},\bv{k})\right]
\nonumber \\
&\hspace{-16cm}+&\hspace{-8cm}
\frac{1}{i\hbar\sqrt{\varepsilon}}\int_{\mathbb{R}^d}\rdf{}{\bv{k}'}{(2\pi)^d}
\re{i\bv{k}'\cdot\bv{x}/\varepsilon}\tilde{V}(\bv{k}')
\left[W_{\varepsilon}(t,\bv{x},\bv{k}-\bv{k}')-
W_{\varepsilon}(t,\bv{x},\bv{k})\right].
\label{eqWeps}
\eea
We note that the term $-\nabla_{\bv{x}}H$ is responsible for the 
Lorentz force:
\bea
-\nabla_{\bv{x}}H
&=&
-\frac{e}{m_{\rm e}}\nabla_{\bv{x}}(\bv{k}\cdot\bv{A})+O(|\bv{A}|^2),
\nonumber \\
-\frac{e}{m_{\rm e}}\nabla_{\bv{x}}(\bv{k}\cdot\bv{A})
&=&
-\frac{e}{m_{\rm e}}\left[
(\bv{k}\cdot\nabla_{\bv{x}})\bv{A}+
\bv{k}\times(\nabla_{\bv{x}}\times\bv{A})\right]
\nonumber \\
&=&
-\frac{e}{\hbar}\left[-\bv{E}+\bv{v}\times\bv{B}\right],
\eea
where $\hbar\bv{k}=m_{\rm e}\bv{v}+e\bv{A}$ and $\bv{v}=\rd{\bv{x}}/\rd{t}$.  
The left-hand side of \eref{eqWeps} can be expressed as
\be
\frac{\partial}{\partial t}W_{\varepsilon}+\left[
H\left(
\bv{x}-\frac{1}{2}\nabla_{\bv{\mu}},\bv{\mu}+\frac{1}{2}\nabla_{\bv{x}},\bv{q}
\right)
-
H\left(
\bv{x}+\frac{1}{2}\nabla_{\bv{\mu}},\bv{\mu}-\frac{1}{2}\nabla_{\bv{x}},\bv{q}
\right)
\right]W_{\varepsilon}(t,\bv{x},\bv{k}).
\ee
This expression was first obtained by Kubo \cite{Kubo64}.

We introduce a two-scale expansion for $W_{\varepsilon}$:
\be
W_{\varepsilon}(t,\bv{x},\bv{k})=
W_0(t,\bv{x},\bv{z},\bv{k})+
\sqrt{\varepsilon}W_1(t,\bv{x},\bv{z},\bv{k})+
\varepsilon W_2(t,\bv{x},\bv{z},\bv{k})+\cdots.
\ee
We assume that $W_0$ is deterministic and periodic with respect to 
$\bv{z}=\bv{x}/\varepsilon$.  We replace
\be
\nabla_{\bv{x}}\to\nabla_{\bv{x}}+\frac{1}{\varepsilon}\nabla_{\bv{z}}.
\ee
By collecting terms of $O(\varepsilon^{-1})$, we have
\be
\mathcal{L}W_0=0,
\label{order-1}
\ee
where the skew symmetric operator $\mathcal{L}$ is given by
\bea
\mathcal{L}f(\bv{z},\bv{k})
&=&
\nabla_{\bv{\mu}}H(\bv{x},\bv{\mu},\bv{q})\cdot\nabla_{\bv{z}}
f(\bv{z},\bv{k})+
\frac{i\hbar}{2m_{\rm e}}\Delta_{\bv{z}}f(\bv{z},\bv{k})
\nonumber \\
&-&
\frac{1}{i\hbar}\sum_{\bv{\mu}'\in L^*}\re{i\bv{\mu}'\cdot\bv{z}}
\tilde{U}(\bv{\mu}')\left[f(\bv{z},\bv{k}-\bv{\mu}')-f(\bv{z},\bv{k})\right]
\nonumber \\
&=&
\left[\frac{\hbar\bv{k}}{m_{\rm e}}+\frac{e}{m_{\rm e}}\bv{A}(t,\bv{x})\right]
\cdot\nabla_{\bv{z}}f(\bv{z},\bv{k})+
\frac{i\hbar}{2m_{\rm e}}\Delta_{\bv{z}}f(\bv{z},\bv{k})
\nonumber \\
&-&
\frac{1}{i\hbar}\sum_{\bv{\mu}'\in L^*}\re{i\bv{\mu}'\cdot\bv{z}}
\tilde{U}(\bv{\mu}')\left[f(\bv{z},\bv{k}-\bv{\mu}')-f(\bv{z},\bv{k})\right].
\eea
By collecting terms of $O(\varepsilon^{-1/2})$, we have
\be
\mathcal{L}W_1=
\frac{1}{i\hbar}\int_{\mathbb{R}^d}\rdf{}{\bv{k}'}{(2\pi)^d}
\re{i\bv{k}'\cdot\bv{z}}\tilde{V}(\bv{k}')
\left[W_0(t,\bv{x},\bv{k}-\bv{k}')-W_0(t,\bv{x},\bv{k})\right].
\label{order-0.5}
\ee
By collecting terms of $O(1)$, we have
\bea
\frac{\partial}{\partial t}W_0+\left(\nabla_{\bv{\mu}}H\cdot\nabla_{\bv{x}}-
\nabla_{\bv{x}}H\cdot\nabla_{\bv{\mu}}\right)W_0
&=&
-\mathcal{L}W_2-
\frac{i\hbar}{m_{\rm e}}\nabla_{\bv{x}}\cdot\nabla_{\bv{z}}W_0(t,\bv{x},\bv{k})
\nonumber \\
&\hspace{-6cm}+&\hspace{-3cm}
\frac{1}{i\hbar}\int_{\mathbb{R}^d}\rdf{}{\bv{k}'}{(2\pi)^d}
\re{i\bv{k}'\cdot\bv{z}}\tilde{V}(\bv{k}')
\left[W_1(t,\bv{x},\bv{k}-\bv{k}')-W_1(t,\bv{x},\bv{k})\right].
\label{order0}
\eea

\section{The Bloch Functions}
\label{sec:bloch}

To obtain the eigenvalues and eigenfunctions of the operator $\mathcal{L}$, 
we consider the following eigenproblem.
\be
\left[\frac{1}{2m_{\rm e}}\left(\frac{\hbar}{i}\nabla_{\bv{z}}+
e\bv{A}(t,\bv{x})\right)^2+U(\bv{z})\right]\Phi_{m\alpha}(\bv{z},\bv{q})=
E_m(\bv{q})\Phi_{m\alpha}(\bv{z},\bv{q}).
\ee
Here, $\alpha$ labels degenerate energy levels of $E_m$ with 
multiplicity $r_m$: $\alpha=1,\dots,r_m$.  We assume that 
there is no level crossing and 
\be
E_1(\bv{q})\le E_2(\bv{q})\le\cdots.
\ee
The parameter $\bv{q}\in\mathbb{R}^d$ labels eigenvalues of 
the translational operator 
$T$ ($=\re{\bv{\nu}\cdot\nabla_{\bv{z}}}$, $\bv{\nu}\in L$):
\be
T\Phi(\bv{z},\bv{q})=\Phi(\bv{z}+\bv{\nu},\bv{q})=
\re{i\bv{q}\cdot\bv{\nu}}\Phi(\bv{z},\bv{q}).
\label{transl}
\ee
Note that $\bv{q}$ moves inside the Brillouin zone ($\Omega_{\rm BZ}$).  
It is, however, convenient to extend $\Phi(\bv{z},\bv{q})$ to $\mathbb{R}^d$ 
in $\bv{q}$ with $L^*$-periodic.
The eigenfunctions $\Phi_{m\alpha}(\bv{z},\bv{q})$ form a complete 
orthonormal basis in $L^2(C)$:
\be
(\Phi_{m\alpha},\Phi_{j\beta})=\int_C\rdf{}{\bv{z}}{|C|}
\Phi_{m\alpha}(\bv{z},\bv{q})\bar{\Phi}_{j\beta}(\bv{z},\bv{q})=
\delta_{mj}\delta_{\alpha\beta}.
\label{Phiortho0}
\ee
We rewrite the eigenproblem using the periodic function
\be
\phi(\bv{z},\bv{q})=\re{-i\bv{q}\cdot\bv{z}}\Phi(\bv{z},\bv{q}).
\ee
We obtain
\be
\left[\frac{1}{2m_{\rm e}}\left(\frac{\hbar}{i}\nabla_{\bv{z}}+
\hbar\bv{q}+e\bv{A}(t,\bv{x})\right)^2+U(\bv{z})\right]
\phi_{m\alpha}(\bv{z},\bv{q})=E_m(\bv{q})\phi_{m\alpha}(\bv{z},\bv{q}).
\ee
By differentiating the equation with respect to $q_j$, we obtain
\be
\frac{\partial E_m}{\partial q_j}\delta_{mn}\delta_{\alpha\beta}=
-i\frac{\hbar^2}{m_{\rm e}}
\left(\frac{\partial\phi_{m\alpha}}{\partial z_j},\phi_{n\beta}\right)+
\frac{\hbar^2}{m_{\rm e}}\left(q_j+\frac{e}{\hbar}A_j\right)
\delta_{mn}\delta_{\alpha\beta}.
\ee
Thus we have
\bea
\frac{\partial E_m}{\partial q_j}\delta_{mn}\delta_{\alpha\beta}
&=&
-i\frac{\hbar^2}{m_{\rm e}}
\left(\frac{\partial\Phi_{m\alpha}}{\partial z_j},\Phi_{n\beta}\right)+
\frac{e\hbar}{m_{\rm e}}A_j(t,\bv{x})\delta_{mn}\delta_{\alpha\beta}.
\eea
Similarly, we also have
\be
\left(A_l\frac{\partial\Phi_{m\alpha}}{\partial z_j},\Phi_{n\beta}\right)=
\frac{im_{\rm e}}{\hbar^2}\frac{\partial E_m}{\partial q_j}
\left(A_l\Phi_{m\alpha},\Phi_{n\beta}\right)-
\frac{ie}{\hbar}A_j(t,\bv{x})\left(A_l\Phi_{m\alpha},\Phi_{n\beta}\right),
\ee
where
\be
\left(A_l\Phi_{m\alpha},\Phi_{n\beta}\right)=
\int_C\rdf{}{\bv{z}}{|C|}A_l(t,\bv{z})
\Phi_{m\alpha}(\bv{z},\bv{q})\bar{\Phi}_{n\beta}(\bv{z},\bv{q}).
\ee
The Bloch functions satisfy the following orthogonality relations 
\cite{Bal99}.
\bea
\frac{1}{\Omega_{\rm BZ}}\sum_{m,\alpha}\int_{\rm BZ}\rd{\bv{q}}
\Phi_{m\alpha}(\bv{x},\bv{q})\bar{\Phi}_{m\alpha}(\bv{y},\bv{q})
&=&
\delta(\bv{y}-\bv{x}),
\label{Phiortho1}
\\
\frac{1}{\Omega_{\rm BZ}}\int_{\mathbb{R}^d}\rd{\bv{x}}
\Phi_{j\alpha}(\bv{x},\bv{q})\bar{\Phi}_{m\beta}(\bv{x},\bv{q}')
&=&
\delta_{jm}\delta_{\alpha\beta}\delta(\bv{q}-\bv{q}').
\label{Phiortho2}
\eea

\section{Decomposition of $W_0$}
\label{sec:w0decomp}

Let us define the $\bv{z}$-periodic functions 
$Q_{mn}^{\alpha\beta}(\bv{z},\bv{\mu},\bv{q})$, $\bv{\mu}\in L^*$, 
$\bv{q}\in\Omega_{\rm BZ}$ by
\be
Q_{mn}^{\alpha\beta}(\bv{z},\bv{\mu},\bv{q})=\int_C\rdf{}{\bv{y}}{|C|}
\re{i(\bv{q}+\bv{\mu})\cdot\bv{y}}
\Phi_{m\alpha}\left(\bv{z}-\bv{y},\bv{q}\right)
\bar{\Phi}_{n\beta}\left(\bv{z},\bv{q}\right).
\ee
These functions satisfy the following orthogonality relation.
\be
\sum_{\bv{\mu}\in L^*}\int_C\rdf{}{\bv{z}}{|C|}
\bar{Q}_{j}^{\alpha\beta}(\bv{z},\bv{\mu},\bv{q})
Q_{m}^{\alpha'\beta'}(\bv{z},\bv{\mu},\bv{q})
=
\delta_{jm}\delta_{\alpha\alpha'}\delta_{\beta\beta'}.
\ee
They are eigenfunctions of $\mathcal{L}$:
\be
\mathcal{L}Q_{mn}^{\alpha\beta}(\bv{z},\bv{\mu},\bv{q})=
\frac{i}{\hbar}\left[E_m(\bv{q})-E_n(\bv{q})\right]
Q_{mn}^{\alpha\beta}(\bv{z},\bv{\mu},\bv{q}),
\label{eigenL}
\ee
with $\bv{k}=\bv{q}+\bv{\mu}$.  

Let us write $Q_{mm}^{\alpha\beta}$ as $Q_{m}^{\alpha\beta}$.  
\Eref{eigenL} implies that $\ker\mathcal{L}$ is spanned by 
$Q_{m}^{\alpha\beta}$.  By \eref{order-1}, $W_0(t,\bv{x},\bv{z},\bv{k})$ 
may be decomposed as
\be
W_0(t,\bv{x},\bv{z},\bv{k})=W_0(t,\bv{x},\bv{z},\bv{q}+\bv{\mu})=
\sum_{m,\alpha,\beta}\{u_m(t,\bv{x},\bv{q})\}_{\alpha\beta}
Q_{m}^{\alpha\beta}(\bv{z},\bv{\mu},\bv{q}),
\ee
where $u_m(t,\bv{x},\bv{q})$ is a $r_m\times r_m$ matrix.

\section{Decomposition of $W_1$}
\label{sec:w1decomp}

Let us look at \eref{order-0.5}.  In general, $W_1$ is not periodic in 
the fast variable $\bv{z}$.  Hence, instead of $Q_{mn}^{\alpha\beta}$, 
we use the basis functions
\be
P_{mn}^{\alpha\beta}(\bv{z},\bv{\mu},\bv{q},\bv{q}_0)=
\int_C\rdf{}{\bv{y}}{|C|}\re{i(\bv{q}+\bv{\mu})\cdot\bv{y}}
\Phi_{m\alpha}\left(\bv{z}-\bv{y},\bv{q}\right)
\bar{\Phi}_{n\beta}\left(\bv{z},\bv{q}+\bv{q}_0\right),
\ee
where $\bv{z}\in\mathbb{R}^d$ and $\bv{q},\bv{q}_0\in\Omega_{\rm BZ}$.  
The functions $P_{mn}^{\alpha\beta}$ are quasi-periodic in $\bv{z}$:
\be
P_{mn}^{\alpha\beta}(\bv{z}+\bv{\nu},\bv{\mu},\bv{q},\bv{q}_0)=
P_{mn}^{\alpha\beta}(\bv{z},\bv{\mu},\bv{q},\bv{q}_0)
\re{-i\bv{\nu}\cdot\bv{q}_0},
\ee
where $\bv{\nu}\in L$.  
Using \eref{musum}, \eref{Phiortho0}, and \eref{Phiortho2}, we obtain 
the following orthogonality relation.
\be
\sum_{\mu\in L^*}\int_{\mathbb{R}^d}\rdf{}{\bv{z}}{\Omega_{\rm BZ}}
P_{mn}^{\alpha\beta}(\bv{z},\bv{\mu},\bv{q},\bv{q}')
\bar{P}_{jl}^{\alpha'\beta'}(\bv{z},\bv{\mu},\bv{q},\bv{q}'')
=
\delta_{mj}\delta_{nl}\delta_{\alpha\alpha'}\delta_{\beta\beta'}
\delta(\bv{q}'-\bv{q}'').
\ee
The functions $P_{mn}^{\alpha\beta}$ are also eigenfunctions of 
$\mathcal{L}$.
\be
\mathcal{L}P_{mn}^{\alpha\beta}(\bv{z},\bv{\mu},\bv{q},\bv{q}_0)
=
\frac{i}{\hbar}\left[E_m(\bv{q})-E_n(\bv{q}+\bv{q}_0)\right]
P_{mn}^{\alpha\beta}(\bv{z},\bv{\mu},\bv{q},\bv{q}_0).
\ee
We write $W_1$ in this basis as
\be
W_1(t,\bv{x},\bv{z},\bv{q}+\bv{\mu})=\sum_{mn\alpha\beta}
\int_{\rm BZ}\rdf{}{\bv{q}'}{\Omega_{\rm BZ}}
\eta_{mn}^{\alpha\beta}(t,\bv{x},\bv{q},\bv{q}')
P_{mn}^{\alpha\beta}(\bv{z},\bv{\mu},\bv{q},\bv{q}'),
\label{W1P}
\ee
where $\bv{z}\in\mathbb{R}^d$, $\bv{q}\in\Omega_{\rm BZ}$, and 
$\bv{\mu}\in L^*$.  We plug this into \eref{order-0.5}, multiply 
$\bar{P}_{jl}^{\alpha\beta}(\bv{z},\bv{\mu},\bv{q},\bv{q}_0)$, 
sum over $\bv{\mu}\in L^*$, and integrate over $\bv{z}\in\mathbb{R}^d$.
The lhs becomes
\bea
\sum_{\bv{\mu}\in L^*}\int_{\mathbb{R}^d}\rd{\bv{z}}
\bar{P}_{jl}^{\alpha\beta}(\bv{z},\bv{\mu},\bv{q},\bv{q}_0)
\mathcal{L}\sum_{mn\alpha'\beta'}\int_{\rm BZ}\rdf{}{\bv{q}'}{\Omega_{\rm BZ}}
\eta_{mn}^{\alpha'\beta'}(t,\bv{x},\bv{q},\bv{q}')
P_{mn}^{\alpha'\beta'}(\bv{z},\bv{\mu},\bv{q},\bv{q}')
&&
\nonumber \\
&\hspace{-14cm}=&\hspace{-7cm}
\frac{i}{\hbar}\left[E_j(\bv{q})-E_l(\bv{q}+\bv{q}_0)\right]
\eta_{jl}^{\alpha\beta}(t,\bv{x},\bv{q},\bv{q}_0).
\eea
We define \cite{Bal99}
\be
T_{jm}^{\alpha\beta}(\bv{q}',\bv{q})=
\int_C\rdf{}{\bv{y}}{(2\pi)^{(d-1)/2}|C|}\re{i(\bv{q}'-\bv{q})\cdot\bv{y}}
\Phi_{m\beta}(\bv{y},\bv{q})\bar{\Phi}_{j\alpha}(\bv{y},\bv{q}').
\ee
By taking the sum over $\bv{\mu}$, the rhs becomes
\bea
\frac{1}{i\hbar}\sum_{\bv{\mu}'\in L^*}\int_{\rm BZ}\rdf{}{\bv{q}'}{(2\pi)^d}
\tilde{V}(\bv{q}'+\bv{\mu}')
\sum_{m\alpha'\beta'}\int_{\mathbb{R}^d}\rd{\bv{z}}
\Biggl[
&&
\nonumber \\
&\hspace{-12cm}&\hspace{-6cm}
(2\pi)^{(d-1)/2}
T_{jm}^{\alpha\alpha'}(\bv{q},\bv{q}-\bv{q}'-\bv{\mu}')
\{u_m(t,\bv{x},\bv{q}-\bv{q}')\}_{\alpha'\beta'}
\bar{\Phi}_{m\beta'}(\bv{z},\bv{q}-\bv{q}')
\nonumber \\
&\hspace{-12cm}-&\hspace{-6cm}
\delta_{jm}\delta_{\alpha\alpha'}
\{u_m(t,\bv{x},\bv{q})\}_{\alpha'\beta'}\re{i(\bv{q}'+\bv{\mu}')\cdot\bv{z}}
\bar{\Phi}_{m\beta'}(\bv{z},\bv{q})
\Biggr]\Phi_{l\beta}(\bv{z},\bv{q}+\bv{q}_0)
\nonumber \\
&\hspace{-12cm}=&\hspace{-6cm}
\frac{1}{i\hbar}\sum_{\bv{\mu}\in L^*}\frac{\Omega_{\rm BZ}}{(2\pi)^{(d+1)/2}}
\tilde{V}(-\bv{q}_0-\bv{\mu})
\sum_{\alpha'}T_{jl}^{\alpha\alpha'}(\bv{q},\bv{q}+\bv{q}_0+\bv{\mu})
\{u_l(t,\bv{x},\bv{q}+\bv{q}_0)\}_{\alpha'\beta}
\nonumber \\
&\hspace{-12cm}-&\hspace{-6cm}
\frac{1}{i\hbar}\int_{\mathbb{R}^d}\rdf{}{\bv{k}'}{(2\pi)^d}
\tilde{V}(\bv{k}')\sum_{\beta'}\{u_j(t,\bv{x},\bv{q})\}_{\alpha\beta'}
\int_{\mathbb{R}^d}\rd{\bv{z}}\re{i\bv{k}'\cdot\bv{z}}
\Phi_{l\beta}(\bv{z},\bv{q}+\bv{q}_0)\bar{\Phi}_{j\beta'}(\bv{z},\bv{q})
\Biggr].
\nonumber \\
\eea
Therefore we obtain
\bea
\eta_{jl}^{\alpha\beta}(t,\bv{x},\bv{q},\bv{q}_0)
&=&
\frac{\Omega_{\rm BZ}}{(2\pi)^{(d+1)/2}}\sum_{\bv{\mu}\in L^*}
\frac{\tilde{V}(-\bv{q}_0-\bv{\mu})}
{E_l(\bv{q}+\bv{q}_0)-E_j(\bv{q})+i\xi}
\nonumber \\
&\hspace{-6cm}\times&\hspace{-3cm}
\sum_{\alpha'}\{u_l(t,\bv{x},\bv{q}+\bv{q}_0)\}_{\alpha'\beta}
T_{jl}^{\alpha\alpha'}(\bv{q},\bv{q}+\bv{q}_0+\bv{\mu})
\nonumber \\
&\hspace{-6cm}-&\hspace{-3cm}
\int_{\mathbb{R}^{2d}}\frac{\rd{\bv{z}}\rd{\bv{k}'}}{(2\pi)^d}
\re{i\bv{k}'\cdot\bv{z}}
\frac{\tilde{V}(\bv{k}')\sum_{\beta'}\{u_j(t,\bv{x},\bv{q})\}_{\alpha\beta'}
\Phi_{l\beta}(\bv{z},\bv{q}+\bv{q}_0)\bar{\Phi}_{j\beta'}(\bv{z},\bv{q})}
{E_l(\bv{q}+\bv{q}_0)-E_j(\bv{q})+i\xi},
\label{eta}
\eea
where $\xi$ ($>0$) is infinitesimally small.

\section{Time Evolution of $W_0$}
\label{sec:w0}

Next let us look at \eref{order0}.  We multiply 
$\bar{Q}_{j}^{\alpha\beta}(\bv{z},\bv{\mu},\bv{q})$, integrate both sides 
over $\bv{z}$, and sum over $\bv{\mu}$.
\bea
\sum_{\bv{\mu}\in L^*}\int_C\rdf{}{\bv{z}}{|C|}
\bar{Q}_{j}^{\alpha\beta}(\bv{z},\bv{\mu},\bv{q})
\nonumber \\
&\hspace{-6cm}\times&\hspace{-3cm}
\Biggl[\frac{\partial}{\partial t}W_0+\frac{\hbar}{m_{\rm e}}
\left(\bv{k}+\frac{e}{\hbar}\bv{A}(t,\bv{x})\right)\cdot\nabla_{\bv{x}}W_0-
\frac{e}{m_{\rm e}}\sum_{ln}\left(k_l+\frac{e}{\hbar}A_l(t,\bv{x})\right)
\frac{\partial A_l(t,\bv{x})}{\partial x_n}\partial_{\mu_n}W_0
\Biggr]
\nonumber \\
&\hspace{-6cm}=&\hspace{-3cm}
\sum_{\bv{\mu}\in L^*}\int_C\rdf{}{\bv{z}}{|C|}
\bar{Q}_{j}^{\alpha\beta}(\bv{z},\bv{\mu},\bv{q})
\nonumber \\
&\hspace{-6cm}\times&\hspace{-3cm}
\left[\frac{\partial}{\partial t}W_0+
\sum_l\left(k_l+\frac{e}{\hbar}A_l(t,\bv{x})\right)
\left(\frac{\hbar}{m_{\rm e}}\frac{\partial}{\partial x_l}-
\frac{e}{m_{\rm e}}A_l(t,\nabla_{\bv{\mu}})\right)W_0\right]
\nonumber \\
&\hspace{-6cm}=&\hspace{-3cm}
-\sum_{\bv{\mu}\in L^*}\int_C\rdf{}{\bv{z}}{|C|}
\bar{Q}_{j}^{\alpha\beta}(\bv{z},\bv{\mu},\bv{q})
\left[\mathcal{L}W_2+
\frac{i\hbar}{m_{\rm e}}\nabla_{\bv{x}}\cdot\nabla_{\bv{z}}W_0\right]
\nonumber \\
&\hspace{-6cm}+&\hspace{-3cm}
\sum_{\bv{\mu}\in L^*}\int_C\rdf{}{\bv{z}}{|C|}
\bar{Q}_{j}^{\alpha\beta}(\bv{z},\bv{\mu},\bv{q})
\nonumber \\
&\hspace{-6cm}\times&\hspace{-3cm}
\frac{1}{i\hbar}\int_{\mathbb{R}^d}\rdf{}{\bv{k}'}{(2\pi)^d}
\re{i\bv{k}'\cdot\bv{z}}\tilde{V}(\bv{k}')
\left[W_1(t,\bv{x},\bv{z},\bv{k}-\bv{k}')-W_1(t,\bv{x},\bv{z},\bv{k})\right].
\label{order0integrated}
\eea
The first integral on the rhs vanishes ($\mathcal{L}$ is skew symmetric and 
$Q_{j}^{\alpha\beta}\in\ker\mathcal{L}$).

\subsection{LHS}
\label{lhs}

The lhs of \eref{order0integrated} is calculated as follows.  
The first term is easy. 
\bea
\sum_{\bv{\mu}\in L^*}\int_C\rdf{}{\bv{z}}{|C|}
\bar{Q}_{j}^{\alpha\beta}(\bv{z},\bv{\mu},\bv{q})
\frac{\partial}{\partial t}W_0
&&
\nonumber \\
&\hspace{-6cm}=&\hspace{-3cm}
\sum_{m,\alpha',\beta'}\frac{\partial}{\partial t}
\{u_m(t,\bv{x},\bv{q})\}_{\alpha'\beta'}
\sum_{\bv{\mu}\in L^*}\int_C\rdf{}{\bv{z}}{|C|}
\bar{Q}_{j}^{\alpha\beta}(\bv{z},\bv{\mu},\bv{q})
Q_{m}^{\alpha'\beta'}(\bv{z},\bv{\mu},\bv{q})
\nonumber \\
&\hspace{-6cm}=&\hspace{-3cm}
\frac{\partial}{\partial t}\{u_j(t,\bv{x},\bv{q})\}_{\alpha\beta}.
\eea
The second term is calculated as follows.
\bea
\sum_{\bv{\mu}\in L^*}\int_C\rdf{}{\bv{z}}{|C|}
\bar{Q}_{j}^{\alpha\beta}(\bv{z},\bv{\mu},\bv{q})
\left(k_l+\frac{e}{\hbar}A_l(\bv{x})\right)
\left(\frac{\partial}{\partial x_l}-\frac{e}{\hbar}
A_l(\nabla_{\bv{\mu}})\right)
\nonumber \\
&\hspace{-16cm}\times&\hspace{-8cm}
\{u_m(t,\bv{x},\bv{q})\}_{\alpha'\beta'}
Q_{m}^{\alpha'\beta'}(\bv{z},\bv{\mu},\bv{q})
\nonumber \\
&\hspace{-16cm}=&\hspace{-8cm}
\sum_{\bv{\mu}\in L^*}\int_C\rdf{}{\bv{z}}{|C|}
\bar{Q}_{j}^{\alpha\beta}(\bv{z},\bv{\mu},\bv{q})
\int_C\rdf{}{\bv{y}'}{|C|}
\left[\left(-i\frac{\partial}{\partial y'_l}+\frac{e}{\hbar}A_l(\bv{x})\right)
\re{i(\bv{q}+\bv{\mu})\cdot\bv{y}'}\right]
\nonumber \\
&\hspace{-16cm}\times&\hspace{-8cm}
\left\{\frac{\partial u_m}{\partial x_l}-
\frac{ie}{\hbar}u_mA_l(\bv{y}')\right\}_{\alpha'\beta'}
\Phi_{m\alpha'}(\bv{z}-\bv{y}',\bv{q})\bar{\Phi}_{m\beta'}(\bv{z},\bv{q}).
\nonumber \\
&\hspace{-16cm}=&\hspace{-8cm}
\frac{\partial \{u_j\}_{\alpha'\beta'}}{\partial x_l}\delta_{jm}
\delta_{\beta\beta'}
\left[
-i\left(\frac{\partial\Phi_{j\alpha'}}{\partial z_l},\Phi_{j\alpha}\right)
+\frac{e}{\hbar}A_l(\bv{x})\delta_{\alpha\alpha'}
\right]
\nonumber \\
&\hspace{-16cm}-&\hspace{-8cm}
\frac{ie}{\hbar}\{u_m\}_{\alpha'\beta'}
\Biggl[
-i\left(\frac{\partial\Phi_{m\alpha'}}{\partial z_l},\Phi_{j\alpha}\right)
\left(A_l\Phi_{j\beta},\Phi_{m\beta'}\right)
+
\frac{e}{\hbar}A_l(\bv{x})\delta_{jm}\delta_{\alpha\alpha'}
\left(A_l\Phi_{j\beta},\Phi_{m\beta'}\right)
\nonumber \\
&\hspace{-16cm}+&\hspace{-8cm}
i\delta_{jm}\delta_{\beta\beta'}
\left(A_l\frac{\partial\Phi_{m\alpha'}}{\partial z_l},\Phi_{j\alpha}\right)
-
\delta_{jm}\delta_{\beta\beta'}\frac{e}{\hbar}A_l(\bv{x})
\left(A_l\Phi_{m\alpha'},\Phi_{j\alpha}\right)
\Biggr].
\eea
Therefore, the lhs of \eref{order0integrated} is written as
\bea
\frac{\partial}{\partial t}\{u_j(t,\bv{x},\bv{q})\}_{\alpha\beta}
+\frac{1}{\hbar}\nabla_{\bv{q}}E_j(\bv{q})\cdot\nabla_{\bv{x}}
\{u_j(t,\bv{x},\bv{q})\}_{\alpha\beta}+\frac{e}{\hbar}
\nonumber \\
&\hspace{-10cm}\times&\hspace{-5cm}
\Biggl[\sum_{\alpha'}\left(
\frac{i}{\hbar}\nabla_{\bv{q}}E_j\cdot\bv{A}(\bv{z})\Phi_{j\alpha'}(\bv{z}),
\Phi_{j\alpha}(\bv{z})\right)\{u_j(t,\bv{x},\bv{q})\}_{\alpha'\beta}
\nonumber \\
&\hspace{-10cm}-&\hspace{-5cm}
\sum_{\beta'}\{u_j(t,\bv{x},\bv{q})\}_{\alpha\beta'}\left(
\frac{i}{\hbar}\nabla_{\bv{q}}E_j\cdot\bv{A}(\bv{z})\Phi_{j\beta}(\bv{z}),
\Phi_{j\beta'}(\bv{z})\right)\Biggr].
\label{order0integratedLHS}
\eea

\subsection{RHS}
\label{sec:rhs}

Let us calculate the rhs of \eref{order0integrated}.  After taking 
ensemble average, the rhs is given by the sum of two matrices 
$I_1$ and $I_2$ where
\bea
\{I_1\}_{\alpha\beta}
&=&
\frac{1}{i\hbar}\sum_{\bv{\mu}\in L^*}\int_C\rdf{}{\bv{z}}{|C|}
\bar{Q}_{j}^{\alpha\beta}(\bv{z},\bv{\mu},\bv{q})
\int_{\mathbb{R}^d}\rdf{}{\bv{k}'}{(2\pi)^d}\re{i\bv{k}'\cdot\bv{z}}
\nonumber \\
&\times&
\langle\tilde{V}(\bv{k}')W_1(t,\bv{x},\bv{z},\bv{k}-\bv{k}')\rangle,
\\
\{I_2\}_{\alpha\beta}
&=&
-\frac{1}{i\hbar}\sum_{\bv{\mu}\in L^*}\int_C\rdf{}{\bv{z}}{|C|}
\bar{Q}_{j}^{\alpha\beta}(\bv{z},\bv{\mu},\bv{q})
\int_{\mathbb{R}^d}\rdf{}{\bv{k}'}{(2\pi)^d}\re{i\bv{k}'\cdot\bv{z}}
\nonumber \\
&\times&
\langle\tilde{V}(\bv{k}')W_1(t,\bv{x},\bv{z},\bv{k})\rangle.
\eea
Using \eref{W1P} and \eref{eta}, $I_1$ is written as
\be
I_1=I_{11}-I_{12},
\ee
where
\bea
\{I_{11}\}_{\alpha\beta}
&=&
\frac{1}{i\hbar}\sum_{\bv{\mu}\in L^*}\int_C\rdf{}{\bv{z}}{|C|}
\bar{Q}_{j}^{\alpha\beta}(\bv{z},\bv{\mu},\bv{q})
\int_{\mathbb{R}^d}\rdf{}{\bv{k}'}{(2\pi)^d}\re{i\bv{k}'\cdot\bv{z}}
\sum_{mn\alpha'\beta'}\int_{\rm BZ}\rdf{}{\bv{q}''}{\Omega_{\rm BZ}}
\nonumber \\
&\times&
\frac{\Omega_{\rm BZ}}{(2\pi)^{(d+1)/2}}\sum_{\bv{\mu}''\in L^*}
\frac{\langle\tilde{V}(\bv{k}')\tilde{V}(-\bv{q}''-\bv{\mu}'')\rangle}
{E_n(\bv{q}-\bv{q}'+\bv{q}'')-E_m(\bv{q}-\bv{q}')+i\xi}
\nonumber \\
&\times&
\sum_{\alpha''}\{u_n(t,\bv{x},\bv{q}-\bv{q}'+\bv{q}'')\}_{\alpha''\beta'}
T_{mn}^{\alpha'\alpha''}(\bv{q}-\bv{q}',\bv{q}-\bv{q}'+\bv{q}''+\bv{\mu}'')
\nonumber \\
&\times&
P_{mn}^{\alpha'\beta'}(\bv{z},\bv{\mu}-\bv{\mu}',\bv{q}-\bv{q}',\bv{q}''),
\eea
and
\bea
\{I_{12}\}_{\alpha\beta}
&=&
\frac{1}{i\hbar}\sum_{\bv{\mu}\in L^*}\int_C\rdf{}{\bv{z}}{|C|}
\bar{Q}_{j}^{\alpha\beta}(\bv{z},\bv{\mu},\bv{q})
\int_{\mathbb{R}^d}\rdf{}{\bv{k}'}{(2\pi)^d}\re{i\bv{k}'\cdot\bv{z}}
\sum_{mn\alpha'\beta'}\int_{\rm BZ}\rdf{}{\bv{q}''}{\Omega_{\rm BZ}}
\nonumber \\
&\times&
\int_{\mathbb{R}^{2d}}\frac{\rd{\bv{z}'}\rd{\bv{k}_1}}{(2\pi)^d}
\re{i\bv{k}_1\cdot\bv{z}'}
\langle\tilde{V}(\bv{k}')\tilde{V}(\bv{k}_1)\rangle
P_{mn}^{\alpha'\beta'}(\bv{z},\bv{\mu}-\bv{\mu}',\bv{q}-\bv{q}',\bv{q}'')
\nonumber \\
&\times&
\frac{\sum_{\beta''}\{u_m(t,\bv{x},\bv{q}-\bv{q}')\}_{\alpha'\beta''}
\Phi_{n\beta'}(\bv{z}',\bv{q}-\bv{q}'+\bv{q}'')
\bar{\Phi}_{m\beta''}(\bv{z}',\bv{q}-\bv{q}')}
{E_n(\bv{q}-\bv{q}'+\bv{q}'')-E_m(\bv{q}-\bv{q}')+i\xi}.
\nonumber \\
\eea
We obtain
\bea
\{I_{11}\}_{\alpha\beta}
&=&
\frac{1}{i\hbar}\sum_{m}
\sum_{\bv{\mu}''\in L^*}\int_{\rm BZ}\rdf{}{\bv{q}''}{2\pi}
\frac{\tilde{R}(-\bv{q}''-\bv{\mu}'')}
{E_j(\bv{q})-E_m(\bv{q}-\bv{q}'')+i\xi}
\sum_{\alpha'\alpha''}\{u_j(t,\bv{x},\bv{q})\}_{\alpha''\beta}
\nonumber \\
&\times&
T_{mj}^{\alpha'\alpha''}(\bv{q}-\bv{q}'',\bv{q}+\bv{\mu}'')
T_{jm}^{\alpha\alpha'}(\bv{q}+\bv{\mu}'',\bv{q}-\bv{q}'').
\eea
Similarly we also obtain $I_{12}$.
\bea
\{I_{12}\}_{\alpha\beta}
&=&
\frac{1}{i\hbar}\int_{\mathbb{R}^{d}}\rdf{}{\bv{k}'}{2\pi}\sum_{m}
\frac{\tilde{R}(-\bv{k}')}{E_j(\bv{q})-E_m(\bv{q}-\bv{q}')+i\xi}
\nonumber \\
&\times&
\sum_{\alpha'\beta''}\{u_m(t,\bv{x},\bv{q}-\bv{q}')\}_{\alpha'\beta''}
T_{jm}^{\alpha\alpha'}(\bv{q},\bv{q}-\bv{k}')
T_{mj}^{\beta''\beta}(\bv{q}-\bv{k}',\bv{q}),
\eea
where we used
\be
\int_{\mathbb{R}^d}\rd{\bv{z}'}=\sum_{\bv{\nu}\in L}\int_C\rd{\bv{y}},
\quad
\frac{1}{\Omega_{\rm BZ}}\sum_{\bv{\nu}\in L}\re{i\bv{q}\cdot\bv{\nu}}
=\sum_{\bv{\mu}\in L^*}\delta(\bv{q}+\bv{\mu}).
\ee
We have
\bea
\{I_1\}_{\alpha\beta}
&=&
\frac{1}{i\hbar}\sum_{m}
\sum_{\bv{\mu}'\in L^*}\int_{\rm BZ}\rdf{}{\bv{q}'}{2\pi}
\frac{\tilde{R}(-\bv{q}'-\bv{\mu}')}{E_j(\bv{q})-E_m(\bv{q}-\bv{q}')+i\xi}
\Biggl[
\nonumber \\
&\hspace{-2cm}&\hspace{-1cm}
\sum_{\alpha'\alpha''}
T_{jm}^{\alpha\alpha'}(\bv{q},\bv{q}-\bv{q}'-\bv{\mu}')
T_{mj}^{\alpha'\alpha''}(\bv{q}-\bv{q}'-\bv{\mu}',\bv{q})
\{u_j(t,\bv{x},\bv{q})\}_{\alpha''\beta}
\nonumber \\
&\hspace{-2cm}-&\hspace{-1cm}
\sum_{\alpha'\beta''}
T_{jm}^{\alpha\alpha'}(\bv{q},\bv{q}-\bv{q}'-\bv{\mu}')
\{u_m(t,\bv{x},\bv{q}-\bv{q}')\}_{\alpha'\beta''}
T_{mj}^{\beta''\beta}(\bv{q}-\bv{q}'-\bv{\mu}',\bv{q})
\Biggr].
\nonumber \\
\eea
In the same way, we have
\bea
\{I_2\}_{\alpha\beta}
&=&
\frac{-1}{i\hbar}\sum_{n}
\sum_{\bv{\mu}'\in L^*}\int_{\rm BZ}\rdf{}{\bv{q}'}{2\pi}
\frac{\tilde{R}(-\bv{q}'-\bv{\mu}')}
{E_n(\bv{q}+\bv{q}')-E_j(\bv{q})+i\xi}
\Biggl[
\nonumber \\
&&
\sum_{\alpha''\beta'}
T_{jn}^{\alpha\alpha''}(\bv{q},\bv{q}+\bv{q}'+\bv{\mu}')
 \{u_n(t,\bv{x},\bv{q}+\bv{q}')\}_{\alpha''\beta'}
T_{nj}^{\beta'\beta}(\bv{q}+\bv{q}'+\bv{\mu}',\bv{q})
\nonumber \\
&-&
\sum_{\beta'\beta''}\{u_j(t,\bv{x},\bv{q})\}_{\alpha\beta''}
T_{jn}^{\beta''\beta'}(\bv{q},\bv{q}+\bv{q}'+\bv{\mu}')
T_{nj}^{\beta'\beta}(\bv{q}+\bv{q}'+\bv{\mu}',\bv{q})
\Biggr].
\nonumber \\
\eea
Thus we see the relation
\be
I_2=I_1^*,
\ee
where $*$ denotes the Hermitian conjugate.  
By summing $I_1$ and $I_2$, and using the relation
\be
\bar{T}_{mj}^{\beta\alpha}(\bv{q},\bv{q}')=
T_{jm}^{\alpha\beta}(\bv{q}',\bv{q}),
\ee
the rhs of \eref{order0integrated} becomes
\bea
I_1+I_1^*
&=&
\frac{1}{i\hbar}\sum_{m}
\sum_{\bv{\mu}'\in L^*}\int_{\rm BZ}\rdf{}{\bv{q}'}{2\pi}
\tilde{R}(-\bv{q}'-\bv{\mu}')
\Biggl[
\nonumber \\
&\hspace{-2cm}&\hspace{-1cm}
T_{jm}(\bv{q},\bv{q}-\bv{q}'-\bv{\mu}')
u_m(t,\bv{x},\bv{q}-\bv{q}')
T_{mj}(\bv{q}-\bv{q}'-\bv{\mu}',\bv{q})
\nonumber \\
&\hspace{-2cm}\times&\hspace{-1cm}
\left(\frac{1}{E_j(\bv{q})-E_m(\bv{q}-\bv{q}')-i\xi}-
\frac{1}{E_j(\bv{q})-E_m(\bv{q}-\bv{q}')+i\xi}\right)
\nonumber \\
&\hspace{-2cm}+&\hspace{-1cm}
\frac{
T_{jm}(\bv{q},\bv{q}-\bv{q}'-\bv{\mu}')
T_{jm}^*(\bv{q},\bv{q}-\bv{q}'-\bv{\mu}')
u_j(t,\bv{x},\bv{q})
}{E_j(\bv{q})-E_m(\bv{q}-\bv{q}')+i\xi}
\nonumber \\
&\hspace{-2cm}-&\hspace{-1cm}
\frac{
u_j(t,\bv{x},\bv{q})
T_{jm}(\bv{q},\bv{q}-\bv{q}'-\bv{\mu}')
T_{jm}^*(\bv{q},\bv{q}-\bv{q}'-\bv{\mu}')
}{E_j(\bv{q})-E_m(\bv{q}-\bv{q}')-i\xi}
\Biggr].
\eea
We change the variable as $\bv{q}'\to\bv{q}-\bv{q}'$ and obtain
\bea
I_1+I_1^*
&=&
\frac{1}{i\hbar}\sum_{m}
\sum_{\bv{\mu}'\in L^*}\int_{\rm BZ}\rdf{}{\bv{q}'}{2\pi}
\tilde{R}(-\bv{q}'-\bv{\mu}')
\Biggl[
\nonumber \\
&\hspace{-2cm}&\hspace{-1cm}
2\pi i
T_{jm}(\bv{q},\bv{q}'-\bv{\mu}')
u_m(t,\bv{x},\bv{q}')
T_{mj}(\bv{q}'-\bv{\mu}',\bv{q})
\delta\left(E_j(\bv{q})-E_m(\bv{q}')\right)
\nonumber \\
&\hspace{-2cm}+&\hspace{-1cm}
\frac{
T_{jm}(\bv{q},\bv{q}'-\bv{\mu}')
T_{jm}^*(\bv{q},\bv{q}'-\bv{\mu}')
u_j(t,\bv{x},\bv{q})
}{E_j(\bv{q})-E_m(\bv{q}')+i\xi}
\nonumber \\
&\hspace{-2cm}-&\hspace{-1cm}
\frac{
u_j(t,\bv{x},\bv{q})
T_{jm}(\bv{q},\bv{q}'-\bv{\mu}')
T_{jm}^*(\bv{q},\bv{q}'-\bv{\mu}')
}{E_j(\bv{q})-E_m(\bv{q}')-i\xi}
\Biggr].
\label{order0integratedRHS}
\eea

\subsection{Radiative Transport Equation}
\label{sec:rte}

We define a vector $\bv{v}_j(\bv{q})$, a matrix $M_j(\bv{q})$, and 
a superoperator $\hat{\mu}_{{\rm L}}$ as
\bea
\bv{v}_j(\bv{q})
&=&
\frac{1}{\hbar}\nabla_{\bv{q}}E_j(\bv{q}),
\\
\{M_j(\bv{q})\}_{\alpha\beta}
&=&
\left(\frac{ie}{\hbar}\bv{v}_j(\bv{q})\cdot\bv{A}\Phi_{j\beta},
\Phi_{j\alpha}\right),
\\
\hat{\mu}_{{\rm L}}: u_j(\bv{q})
&\mapsto&
[M_j(\bv{q}),u_j(\bv{q})]=M_j(\bv{q})u_j(\bv{q})-u_j(\bv{q})M_j(\bv{q}).
\eea
Furthermore let us define the following superoperators.
\bea
\hat{A}(\bv{q},\bv{q}'): u_j(\bv{q}_0)
&\mapsto&
\frac{1}{\hbar}\sum_m\sum_{\bv{\mu}'\in L^*}\tilde{R}(-\bv{q}'-\bv{\mu}')
\delta\left(E_j(\bv{q})-E_m(\bv{q}')\right)
\nonumber \\
&&
T_{jm}(\bv{q},\bv{q}'-\bv{\mu}')u_m(\bv{q}_0)T_{mj}(\bv{q}'-\bv{\mu}',\bv{q}),
\\
\hat{\mu}_{{\rm s}}: u_j(\bv{q})
&\mapsto&
\frac{-1}{i\hbar}\sum_{m}
\sum_{\bv{\mu}'\in L^*}\int_{\rm BZ}\rdf{}{\bv{q}'}{2\pi}
\tilde{R}(-\bv{q}'-\bv{\mu}')
\nonumber \\
&&
\Biggl[
\frac{
T_{jm}(\bv{q},\bv{q}'-\bv{\mu}')T_{jm}^*(\bv{q},\bv{q}'-\bv{\mu}')u_j(\bv{q})
}{E_j(\bv{q})-E_m(\bv{q}')+i\xi}
\nonumber \\
&-&
\frac{
u_j(\bv{q}_0)T_{jm}(\bv{q},\bv{q}'-\bv{\mu}')T_{jm}^*(\bv{q},\bv{q}'-\bv{\mu}')
}{E_j(\bv{q})-E_m(\bv{q}')-i\xi}
\Biggr].
\eea
Finally,by Eqs.~(\ref{order0}), (\ref{order0integrated}), 
(\ref{order0integratedLHS}), and (\ref{order0integratedRHS}), 
we obtain the RTE as
\bea
\frac{\partial}{\partial t}u_j(t,\bv{x},\bv{q})+
\bv{v}_j(\bv{q})\cdot\nabla_{\bv{x}}u_j(t,\bv{x},\bv{q})+
(\hat{\mu}_{{\rm L}}+\hat{\mu}_{{\rm s}})u_j(t,\bv{x},\bv{q})
&&
\nonumber \\
&&\hspace{-2cm}
=\int_{\rm BZ}\rd{\bv{q}'}
\hat{A}(\bv{q},\bv{q}')u_j(t,\bv{x},\bv{q}').
\label{main}
\eea
Note that the term for $\hat{\mu}_{{\rm L}}$ stems from the 
Lorentz force.  The terms for $\hat{\mu}_{{\rm s}}$ and $\hat{A}$ 
stem from the random potential.




\end{document}